\documentclass[12pt]{article}
\usepackage{amsfonts,amssymb,epsfig}
\usepackage{amsmath}

\date{}
\begin{document}
\newtheorem{df}{Definition}
\newtheorem{thm}{Theorem}

\newtheorem{lm}{Lemma}
\newtheorem{pr}{Proposition}
\newtheorem{co}{Corollary}
\newtheorem{re}{Remark}
\newtheorem{note}{Note}
\newtheorem{claim}{Claim}
\newtheorem{problem}{Problem}

\def\R{{\mathbb R}}

\def\E{\mathbb{E}}
\def\calF{{\cal F}}
\def\N{\mathbb{N}}
\def\calN{{\cal N}}
\def\calH{{\cal H}}
\def\n{\nu}
\def\a{\alpha}
\def\d{\delta}
\def\t{\theta}
\def\e{\varepsilon}
\def\t{\theta}
\def\g{\gamma}
\def\G{\Gamma}
\def\b{\beta}
\def\pf{ \noindent {\bf Proof: \  }}

\newcommand{\qed}{\hfill\vrule height6pt
width6pt depth0pt}
\def\endpf{\qed \medskip} \def\colon{{:}\;}
\setcounter{footnote}{0}

\renewcommand{\qed}{\hfill\vrule height6pt  width6pt depth0pt}

\title{Subspaces of $L_p$ that embed into $L_p(\mu)$ with $\mu$ finite
}

\author{William B. Johnson\thanks{Supported in part by NSF DMS-1001321 and U.S.-Israel Binational Science Foundation
 }  \ and Gideon Schechtman\thanks{Supported in part by U.S.-Israel Binational Science Foundation. Participant NSF Workshop in Analysis and Probability, Texas A\&M University
 } } \maketitle

\begin{abstract}
Enflo and Rosenthal \cite{er} proved that $\ell_p(\aleph_1)$, $1 < p < 2$,  does not (isomorphically) embed into $L_p(\mu)$ with $\mu$ a finite measure.  We prove that if $X$ is a subspace of an $L_p$ space, $1< p < 2$, and $\ell_p(\aleph_1)$ does not embed into $X$, then $X$ embeds into $L_p(\mu)$ for some  finite measure $\mu$.

\end{abstract}

\section{Introduction}
In this note we study the structure of non separable subspaces $X$ of $L_p(\mu)$ with $\mu$ a finite measure.  For $2<p<\infty$ an obvious necessary condition is that $\ell_p(\aleph_1)$  does not (isomorphically) embed into  $X$. Indeed, since every operator from even $\ell_p$ into a Hilbert space is compact, there is no one to one (bounded, linear)  operator from  $\ell_p(\aleph_1)$ into a Hilbert space.  On the other hand, for $2<p$, if $\mu$ is a finite measure we have $L_p(\mu) \subset L_2(\mu)$ with the injection being continuous.  We conjecture that if $X$ is a subspace of some $L_p$ space, $2<p<\infty$, and $\ell_p(\aleph_1)$  does not  embed into  $X$, then $X$ embeds into $L_p(\mu)$ for some  finite measure $\mu$. At the end of Section  3 we verify this conjecture when $X$ is a complemented subspace of  some $L_p$.
However, the only information we know related to this conjecture for general subspaces is
Proposition \ref{p>2}, which says  that in this range of $p$, a subspace $X$ of an  $L_p$ space contains $\ell_p(\aleph_1)$ isomorphically iff $X$ contains $\ell_p(\aleph_1)$ isometrically iff there is no one to one operator from $X$ into a Hilbert space.

For $1\le p < 2$, the conjectured classification mentioned above for $2<p$ is true.  In Theorem \ref{p<2}, the main result of Section 2,  we prove for $p$ in this range that a subspace $X$ of an $L_p$ space embeds into $L_p(\mu)$ for some  finite measure $\mu$ if and only if  $\ell_p(\aleph_1)$  does not embed (isomorphically) into  $X$. This is not equivalent to saying that $X$ does not contain an isometric copy of $\ell_p(\aleph_1)$ (but is equivalent to saying that $X$ contains almost isometric copies of $\ell_p(\aleph_1)$; see the remark after Theorem \ref{p<2}).   Part of Theorem \ref{p<2} is known. The $p=1$ case is an almost immediate consequence of a result due to Rosenthal \cite{ros}, and, as mentioned in the abstract, the fact that $X$ does not embed into  $L_p(\mu)$ for any  finite measure $\mu$ if  $\ell_p(\aleph_1)$ embeds into $X$ is due to Enflo and Rosenthal \cite{er}.  The main new result herein is the ``if" part of Theorem \ref{p<2}.

In Section 3 we generalize the results in Section 2 to higher cardinals and thereby obtain characterizations of subspaces of $L_p$ spaces that contain an isomorphic copy of $\ell_p(\aleph)$ for a general uncountable cardinal $\aleph$. Section 3  contains as special cases the results of Section 2, but the arguments in Section 3 require somewhat more background knowledge and there are a few more technicalities.

Although not directly relevant for this note, we draw attention to another result in \cite{er}; namely, that for $1<p\not= 2 < \infty$, a space $L_p(\mu)$ with $\mu$ finite  does not have an unconditional basis if its density character is at least $\aleph_\omega$. Trying (unfortunately, unsuccessfully) to decide what happens when the density character is $\aleph_1$ led us to the results presented here.

\section{Main Results}
We begin with the easy result mentioned in the introduction

\begin{pr} \label{p>2}
Let $X$ be a subspace of some $L_p$ space, $2<p<\infty$. The following are equivalent:
\begin{itemize}
\item $\ell_p(\aleph_1)$  isometrically  embeds into $X$.
\item There is a subspace of $X$ that is isomorphic to $\ell_p(\aleph_1)$  and is complemented in $L_p$.
\item $\ell_p(\aleph_1)$  isomorphically  embeds into $X$.
\item There is no one to one (bounded, linear) operator from $X$ into a Hilbert space.
\end{itemize}
\end{pr}
\pf
Since every isometric copy of an $L_p$ space in an $L_p$ space is norm one complemented \cite[Theorem 6.3]{lac}, and since we already explained in the introduction why the third assertion implies the fourth,
we only need to prove that the fourth condition implies the first condition.
 For this, we use Maharam's theorem \cite{mah},  \cite[Theorem 5.8]{lac}, which implies that $X$ is a subspace of $L_p := (\sum_{\gamma\in \Gamma} L_p\{-1,1\}^{\aleph_\gamma})_p$ for some set $\Gamma$ of   ordinal numbers, where $\{-1,1\}$ is endowed with the uniform probability measure.

Now assume that there is no one to one operator from $X$ into a Hilbert space.  This implies that for any countable subset $\Gamma'$ of $\Gamma$, the natural projection $P_{\Gamma'}$ from $L_p$ onto $(\sum_{\g \in \Gamma'}  L_p \{-1,1\}^{\aleph_\gamma})_p$
is not one to one on $X$, because one can map $(\sum_{\g \in \Gamma'}  L_p\{-1,1\}^{\aleph_\gamma})_p$ one to one into the Hilbert space
$(\sum_{\g \in \Gamma'}  L_2\{-1,1\}^{\aleph_\gamma})_2$ in an obvious way when  $\Gamma'$ is countable. On the other hand, given any $x$ in $X$, there is a countable subset $x(\Gamma)$ of $\Gamma$ so that $P_\g x=0$ for all $\g$ not in $x(\Gamma)$.  Thus if one takes a collection of unit vectors $x$ in $X$ maximal with respect to the property that $x(\Gamma)\cap y(\Gamma) = \emptyset$ when $x\not= y$, then the collection must have cardinality at least $\aleph_1$ and hence $\ell_p(\aleph_1)$ embeds isometrically into $X$.
\qed

Proposition \ref{p>2} is completely wrong for $p<2$.  For one thing, there is an obvious one to one operator from $\ell_p(\Gamma)$ into a Hilbert space--the formal identity mapping from $\ell_p(\Gamma)$ into $\ell_2(\Gamma)$.  Secondly, there are subspaces of $L_p$ isomorphic to  $\ell_p(\Gamma)$ for any set $\Gamma$ that do not contain isometric copies even of $\ell_p$.  Indeed, take a   family $(f_\g)_{\g\in \Gamma}$ of independent standard normal random variables on some probability space $(\Omega, \mu)$ and in $L_p(\mu)\oplus_p \ell_p(\Gamma)$ consider the closed linear span of $(f_\g \oplus e_\g)_{\g\in \Gamma}$, where $( e_\g)_{\g\in \Gamma}$ is the unit vector basis for $\ell_p(\Gamma)$.

It is reasonable to conjecture that an isomorphic copy of $\ell_p(\Gamma)$ in an $L_p$ space contains for every $\e>0$ a $1+\e$-isomorphic copy of $\ell_p(\Gamma)$, and this is proved in the remark after Theorem \ref{p<2} when $\Gamma$ has cardinality $\aleph_1$, and in Theorem \ref{p<2.1} for general uncountable $\Gamma$.  (The case when  $\Gamma$ is countably infinite is contained in \cite{er}.)

The main result in this note is

\begin{thm} \label{p<2}
Let $X$ be a subspace of some $L_p$ space, $1\le p<2$.  Then $X$ embeds into $L_p(\mu)$ for some  finite measure $\mu$ if and only if  $\ell_p(\aleph_1)$  does not embed (isomorphically) into  $X$.
\end{thm}
\pf
In view of \cite{er}, we only need to prove the if part, so assume that $\ell_p(\aleph_1)$  does not embed  into  $X$.  As in the proof of Proposition \ref{p>2}, by Maharam's theorem we can assume that
$X$ is a subspace of $L_p := (\sum_{\gamma\in \Gamma} L_p\{-1,1\}^{\aleph_\gamma})_p$ for some set $\Gamma$ of ordinal numbers. Assume now that $X$ does not embed
into $L_p(\mu)$ for any  finite measure $\mu$.  We want to build a long unconditionally basic sequence $(x_\a)_{\a<\aleph_1}$ of unit vectors in $X$ that have ``big disjoint pieces"; more precisely, so that there are disjoint countable subsets $\Gamma_\a$, $\a<\aleph_1$,  of $\Gamma$ so that $P_{\Gamma_\a} x_\b \not= 0$ iff $\a = \b$.  It then follows that for some $\e>0$, $\|P_{\Gamma_\a} x_\a\| > \e$ for $\aleph_1$ values of $\a$, which we  might as well assume for all $\a<\aleph_1$.  Such a sequence must, by the diagonal principle \cite[Proposition 1.c.8]{lt}  (or with a worse constant, by a square function argument), dominate the unit vector basis of $\ell_p(\aleph_1)$. But by the type $p$ property of $L_p$ \cite[Theorem 6.2.14]{ak}, every normalized  unconditionally basic sequence of cardinality $\aleph_1$ is dominated by the unit vector basis of $\ell_p(\aleph_1)$.
Here we use ``dominate" as is customary in Banach space theory: $(x_\a)$ dominates $(y_\a)$ provided there is a constant $C$ so that for all finite sets $a_\a$ of scalars,
$\|\sum a_\a y_\a\|\le C \|\sum a_\a x_\a\|   $.

Since  $X$ does not embed
into $L_p(\mu)$ for any  finite measure $\mu$,
we have  for any countable subset $\Gamma'$ of $\Gamma$ that the restriction of $P_{\Gamma'} $ to $X$ is not an isomorphism, because
 $ (\sum_{\gamma\in \Gamma'} L_p\{-1,1\}^{\aleph_\gamma})_p$ is isometrically isomorphic to $L_p(\mu)$ for some finite $\mu$ when $\Gamma'$ is countable. From this it is not hard to get a {\sl set} $(x_\a)_{\a<\aleph_1}$ of unit vectors in $X$ that have big disjoint pieces.  In the case $p=1$, this is enough by Rosenthal's technique \cite{ros}  to get a subset of $(x_\a)_{\a<\aleph_1}$  that is equivalent to the unit vector basis of $\ell_1(\aleph_1)$, but for $p>1$ we need to do more work to get $(x_\a)_{\a<\aleph_1}$ unconditionally basic.  So, from here on, we assume that $1<p<2$.

Call a set $S$ of vectors in $L_p =(\sum_{\gamma\in \Gamma} L_p\{-1,1\}^{\aleph_\gamma})_p$ a {\sl generalized martingale difference set} (GMD set, in short) provided that for every finite subset $F$ of $S$ and every $\g$ in $\G$, the sequence $(P_\g x)_{x\in F}$ can be ordered to be a martingale difference sequence. We allow $0$ to appear in a martingale difference sequence, but the definition requires that $P_\g x \not= P_\g y$ if $P_\g x \not= 0$.  Since a martingale difference sequence is unconditional in $L_p(\mu)$  for any probability $\mu$, any $1<p<\infty$, and with  the unconditional constant depending only on $p$ \cite{bur}, a GMD set in $L_p$ is unconditionally basic for our range of $p$.

Take a collection $V$ of pairs $(x, \g(x))_{x\in M}$ in $X\times \G$ maximal with respect to the properties that $\|x\| =1$, $P_{\g(x)} x \not= 0$, the $\g(x)$ are all distinct, and $M$ is a GMD set.  The collection $M$ is unconditionally basic and has disjoint pieces, hence if $M$ is uncountable there is a subset of $M$ having  cardinality $\aleph_1$ that is equivalent to the unit vector basis for $\ell_p(\aleph_1)$. So assume that $M$ is countable.  For each $x$ in $M$, there is a countable subset $x(\G)$ of $\G$  so that $P_\g x =0$ for $\g\not\in x(\G)$. Set $\G' = \cup_{x\in M} x(\G)$.  Each vector in $L_p\{-1,1\}^{\aleph_\gamma}$ depends on only countably many coordinates, so for each $\g $ in $\G'$  there is countable subset $S(\g)$ of $\aleph_\gamma$ so that for every $x$ in $M$, the vector $P_\g x$ depends only on  $S(\g)$.  Let $Q_S(\g)$ be the (norm one) conditional expectation projection from  $L_p\{-1,1\}^{\aleph_\gamma}$ onto its subspace of functions that depend only on $S(\g)$ and let $Q$ be the (norm one) projection on $L_p$ that is the direct sum over $\G'$ of $Q_S(\g)$. Since the  projection $Q$ has separable range and we are assuming that $X$ does not embed into $L_p(\mu)$ for any finite measure $\mu$, the space $X' := \text{ker} Q\cap X$ also does not  embed into $L_p(\mu)$ for any finite measure $\mu$. This is perhaps not quite obvious but follows from the reflexivity of $X$.  Lindenstrauss proved that reflexive spaces have the separable complementation property \cite{lin}, so there is a complemented separable subspace $Y$ of $X^*$ that contains the range of $(Q_{|X})^*$ and hence the co-separable complemented subspace $Y_\perp\subset X' $ of $X$ cannot embed into  $L_p(\mu)$ for any finite measure $\mu$. Now let $P_{\G'}$ be the natural projection from $L_p$ onto $ (\sum_{\gamma\in \Gamma'} L_p\{-1,1\}^{\aleph_\gamma})_p$. Since $\G'$ is countable, the space  $ (\sum_{\gamma\in \Gamma'} L_p\{-1,1\}^{\aleph_\gamma})_p$ is isometric to $L_p(\mu)$ for some probability $\mu$ and hence $P_{\G'}$ is not an isomorphism on $X'$. So there is a unit vector $x$ in $X'$ so that $P_{\G'}x \not= 0$; in fact, $\|P_{\G'}x\|$ can be taken  arbitrarily close to one (this is important for the remark below). In particular, there is $\g(x)\in \G\sim \G'$ so that $P_{\g(x)} x\not = 0$.  Since also $\{x\}\cup M$ is a GMD set,  $V\cup \{(x,\g(x))\}$ contradicts the maximality of $V$.
\qed

\bigskip

\noindent
{\bf Remark.}
A minor modification of the above argument shows that if $X \subset L_p$, $1<p<2$, does not embed into
$L_p(\mu)$ for some  finite measure $\mu$, then $\ell_p(\aleph_1)$ almost isometrically  embeds   into  $X$.
One needs to define the collection $V$ a bit differently. Fix $\e>0$ and take a set $V$ of pairs $(x,\G(x))_{x\in M}$ with each $x$ a unit vector in $X$, each $\G(x)$ a countable subset of $\G$, and maximal with respect to the properties that the $M$ forms a GMD set,  $\|P_{\G(x)} x\| > 1-\e$, and the $\G(x)$ are pairwise disjoint.  The argument for Theorem \ref{p<2} shows that $V$ is uncountable. Notice that the collection M is an $\e$-GMD perturbation of the disjoint vectors $(P_{\G(x)} x)_{x\in M}$; i.e., $(x -P_{\G(x)} x)_{x\in M}$ is a GMD set of vectors each of which has norm less than $\e$.  Thus for any scalars $a_x$, we have
by the type $p$ property of $L_p$ and the GMD property of $(x -P_{\G(x)} x)_{x\in M}$ that
$\| \sum a_x (x -P_{\G(x)} x) \| \le \e C_p( \sum |a_x|^p)^{1/p}$.

\section{Higher cardinals}
In this section we generalize the main results in the previous section to the setting of $\ell_p(\aleph_\a)$ with $\a >1$. Towards the end of this section we also deal with complemented subspaces
of $L_p$ spaces which do not contain large $ \ell_p(\Gamma)$ spaces and extend
our results to the range $2<p<\infty$, but only for complemented subspaces.
We begin with a higher cardinal version of Proposition \ref{p>2}.

Let us say that a Banach space is an $L_p(\aleph)$ space, where $\aleph$ is an infinite cardinal, provided $X$ is isometric to $(\sum_{\a\in\Gamma} L_p(\mu_\a))_p$ with $|\Gamma|\le \aleph$ and each $\mu_a$ a finite measure (which of course can be taken to be probabilities).

\begin{pr}\label{p>2.1}
Let $X$ be a subspace of some $L_p$ space, $2<p<\infty$, and let $\aleph$ be an uncountable cardinal. The following are equivalent:
\begin{itemize}
\item[(1)]  $\ell_p(\aleph)$  isometrically  embeds into $X$.
\item[(2)] There is a subspace of $X$ that is isomorphic to $\ell_p(\aleph)$  and is complemented in $L_p$.
\item[(3)] $\ell_p(\aleph)$  isomorphically  embeds into $X$.
\item[(4)] There is no one to one (bounded, linear) operator from $X$ into an $L_p(\Gamma)$ space with $\Gamma< \aleph$.
\end{itemize}
\end{pr}

\pf
As was mentioned in the proof of Proposition \ref{p>2}, the implication $(1)\implies (2)$ is known, and $(2)\implies (3)$ is obvious.
For  $(3)\implies (4)$ it is enough to show that for $\G<\aleph$, there is no one to one operator from $\ell_p(\aleph)$  into  $(\sum_{\a\in\Gamma} L_p(\mu_\a))_p$
with each $\mu_\a$ a probability. Suppose, to the contrary, that $T$ is such an operator. Let $i_{p,2}^\a$ be the formal inclusion mapping from $L_p(\mu_\a))_p$ into $L_2(\mu_\a))_p$.  The operators  $i_{p,2}^\a P_\a T$ are all compact and the $i_{p,2}^\a$ are all one to one, so for each $\a\in \G$ there is a countable subset $A_a$ of $\aleph$ such that $P_\a T e_\b=0$ for all $\b \in \aleph\sim A_a$, where $(e_\b)_{\b \in  \aleph}$ is the unit vector basis for $\ell_p(\aleph)$.  Since $T$ is one to one, $\cup_{\a \in \G} = \aleph$ and hence $\aleph = |\cup_{\a \in \G}| \le |\G| \cdot \aleph_0 = |\G|$.

For $(4) \implies (1)$, assume that $X$ is a subspace of $(\sum_{\a\in\Gamma} L_p(\mu_\a))_p$  with each $\mu_\a$ a probability. By (4), for all $\G'\subset \G$ with $|\G'| < \aleph$, we have that the restriction of $P_{\G'}$ to $X$ is not one to one.  Take a collection $\mathcal{S}$ of unit vectors in $X$ maximal with respect to the property that $x(\G)\cap y(\G)=\emptyset$ for $x\not= y$ in $\mathcal{S}$, where $x(\G)$ is the (countable) set of all $\a \in \G$ for which $P_\a x\not= 0$. If $|\mathcal{S}|< \aleph$ then $\G' := \cup_{x \in \mathcal{S}}$ has cardinality at most
$|\mathcal{S}|\cdot \aleph_0 <\aleph$, which by (4) implies that the restriction of $P_{\G'}$ to $X$ is not one to one, which clearly contradicts the maximality of  $\mathcal{S}$.  \qed

We turn now to the case $1<p<2$ (Rosenthal \cite{ros} treated the case $p=1$ long ago). In addition to the argument for Theorem \ref{p<2},
we need the following lemma.

\begin{lm}\label{lemma2}
Let $1<p<2$ and let $\aleph$ be an uncountable cardinal. If $\aleph' < \aleph$, then $\ell_p(\aleph)$ is not isomorphic to a subspace of any $L_p(\aleph')$ space.
\end{lm}

\pf Notice that it follows by duality from Proposition \ref{p>2.1} that $\ell_p(\aleph)$ is not isomorphic to a {\sl complemented} subspace of any $L_p(\aleph')$ space. So we just need to prove that if $\ell_p(\aleph)$ embeds into an $L_p(\aleph')$ space, then it embeds into some other $L_p(\aleph')$ space as a complemented subspace.

Assume that $(x_\a)_{\a\in \aleph}$ is a normalized set of vectors in some $L_p(\aleph')$ space $L_p(\Omega,\mu)$  that is equivalent to the unit vector basis for $\ell_p(\aleph)$.  In particular,  there is $\theta >0$ so that for all finite subsets $F$ of $\aleph$ and all scalars $(c_\a)_{\a \in F}$,
\begin{equation}\label{theta}
\|\sum_{\a\in F} c_\a x_\a\|_p \ge \theta (\sum_{\a\in F} |c_\a|^p)^{1/p}.
\end{equation}

 By enlarging $(\Omega,\mu)$ (but keeping $L_p(\Omega,\mu)$ an $L_p(\aleph')$ space) we can assume that $\mu$ is purely nonatomic.  In the $L_p(\aleph')$ space $L_p(\Omega \times \{-1,1\}^\aleph, \mu \times \nu)$ (where $\nu$ is the usual Haar measure on $\{-1,1\}^\aleph$), consider the vectors $y_\a := x_\a \otimes r_\a$, where $r_\a$ is the usual coordinate projection [Rademacher] on $\{-1,1\}^\aleph$.  Then $(y_\a)$ is also equivalent to the usual basis for $\ell_p(\aleph)$ (and satisfies the inequality (\ref{theta})) and, incidentally, is $1$-unconditional (so that it satisfies the reverse inequality to (\ref{theta}) with $\theta$ replaced by one).

 We claim that the closed linear span $Y$  of $(y_\a)$ is complemented in   $L_p(\Omega \times \{-1,1\}^\aleph, \mu \times \nu)$.  Notice that to prove the claim, it is enough to define for each finite subset $F$ of $\aleph$ a projection $P_F$ from $L_p(\Omega \times \{-1,1\}^\aleph, \mu \times \nu)$ onto the span $Y_F$ of $\{y_\a : \a \in F\}$ so that $\sup_F \|P_F\| : = C <\infty$. Indeed, if you index the finite subsets of $\aleph$ by inclusion, the resulting net $(P_F)$, being uniformly bounded, has a weak operator cluster point (say $P$) by the reflexivity of $L_p$ spaces.  It is easy to check that $P$ is a projection onto $Y$ and $\|P\| \le C$.

So let $F$ be any finite non empty subset of $\aleph$.  By a result of Dor \cite[Theorem B]{dor},  inequality (\ref{theta}) implies that  there are disjoint subsets $(\Omega_\a)_{\a\in F}$
of $\Omega$ so that for each $\a$ in $F$,
\begin{equation}\label{theta1}
\|x_a 1_{\Omega_\a} \|_p\ge \theta^{2/(2-p)}.
\end{equation}
(Dor's theorem is stated for $[0,1]$ with Lebesgue meaure, but the proof works for any non atomic measure space. Alternatively, the more general result can be deduced formally from the case of $[0,1]$.) For simplicity, and without loss of generality by replacing   the $\Omega_\a$ with subsets, we can assume that for each $\a$ in $F$ there is equality in inequality (\ref{theta1}).

Define for $\a$ in $F$ disjoint functions
\begin{equation}\label{fsuba}
f_\a :=       {\theta^{(2p)/(p-2)}}    |x_\a|^{p-1} 1_{\Omega_\a} \text{sign} \, x_\a.
\end{equation}
Now define $g_\a := f_\a \otimes r_\a$.
So the $g_\a$ are disjoint vectors in $L_{p'}(\Omega \times \{-1,1\}^\aleph, \mu \times \nu)$, the dual of $L_p(\Omega \times \{-1,1\}^\aleph, \mu \times \nu)$. The power of $\theta$ in their definition was chosen to make $(y_\a,g_\a)_{\a\in F}$ a biorthogonal system.  A routine computation shows that the projection $P_F := \sum_{\a \in F} g_\a \otimes y_\a$ from  $L_p(\Omega \times \{-1,1\}^\aleph, \mu \times \nu)$ onto $Y_F$ has norm at most   $\theta^{2/(p-2)}$.    \qed

The main result of this section is a generalization of Theorem \ref{p<2} to higher cardinals.

\begin{thm} \label{p<2.1}
Let $X$ be a subspace of some $L_p$ space, $1\le p<2$, and let $\aleph$ be an uncountable cardinal.  The following are equivalent.
\begin{itemize}
\item[(1)]  For all $\e>0$, $\ell_p(\aleph)$ is $1+\e$-isomorphic to a subspace of $X$.
\item[(2)]  There is a subspace of $X$ that is isomorphic to $\ell_p(\aleph)$  and is complemented in $L_p$.
\item[(3)] $\ell_p(\aleph)$  isomorphically  embeds into $X$.
\item[(4)]  $X$ does not isomorphically embed into an $L_p(\aleph')$ space with $\aleph' < \aleph$.
\end{itemize}
\end{thm}
\pf
The implication $(1)\implies (2)$ follows easily from  \cite{sch}, while $(2) \implies (3)$ is obvious.  Lemma \ref{lemma2} gives $(3) \implies (4)$, so we only need to prove $(4) \implies (1)$.

By Maharam's theorem \cite{mah}, we can assume that $X$   is a subspace of $L_p:= (\sum_{\b\in B} L_p\{-1,1\}^{\aleph_\b})_p$.
Given $\e>0$, take a set $M:= \{x_\g : \g\in \G\}$ of unit vectors in $X$ maximal with respect to the properties
\begin{itemize}
\item[(i)] The collection $M$ is a GMD set.
\item[(ii)] There are disjoint countable subsets $B_\g$ of $B$ such that for all $\g$ in $\G$, $\|P_{B_\g} x_\g \| > 1-\e$.
\end{itemize}
To get (1), just as in the proof of Theorem \ref{p<2} it is enough to verify that $|\G|\ge \aleph$. Assume instead that $\aleph > |\G|$. We might as well assume, by enlarging the $B_\g$'s, that for $\b\not\in B' := \cup_{\g\in\G} B_\g$, we have  $P_\b x_\g = 0 $ for all $\g$ in $\G$.

For each $\b$ in $B'$ and $\g$ in $\G$, let $C_{\b,\g}$ be the (countable) set of coordinates of $\{-1,1\}^{\aleph_\b}$ on which $P_\b x_\g$ depends, and let $C_\b := \cup_{\g\in\G} C_{\b,\g}$.  Then $|C_\b|\le |\G|$ and hence the density character of $ L_p\{-1,1\}^{C_\b}$ is at most $|\G|$.

 Let   $\mathbb{E}_\b$  be the (norm one) conditional expectation projection from $L_p\{-1,1\}^{\aleph_\b}$  onto $L_p\{-1,1\}^{C_\b}$ and let $Q:= (\sum_{\b\in B'} \Bbb{E}_\b)_p$ be the direct sum of these projections; we consider $Q$ a a projection on $L_p$. Note that the density character of the range of $Q$ is at most $|\G|<\aleph$.

 Just as in the proof of Theorem \ref{p<2}, this implies that $X':= \ker(Q) \cap X$ does not embed isomorphically into an $L_p(\aleph')$ space with $\aleph'<\aleph$.  Indeed, Lindenstrauss \cite{lin} proved that there is a norm one complemented subspace $Y$ of $X^*$ which contains the range of $(Q_{|X})^*$ so that the density character of $Y$ is the same as the density character of $QX$, which is at most $|\G|<\aleph$. This obviously implies that $Y^*$ (which is isometric to  the range of the adjoint of the norm one projection onto $Y$) embeds isometrically into an $L_p(|\G|)$ space  and hence $Y^\perp$  cannot embed into an $L_p(\aleph')$ space because $X$ does not.  Since  $Y^\perp$  is contained in  $X'$, also $X'$ also
 does not embed into an $L_p(\aleph')$ space.

 The above argument shows that $(P_{B'})_{|X}$ is not an isomorphism, so there is a norm one element $x$ in $X'$ such that $\|P_{B'} x\| < \e$.  Then $M \cup \{x\}$ satisfies (i) and (ii), which contradicts the maximality of $M$.   \qed

  Theorem \ref{p<2.1+} below gives another equivalence to the conditions in Theorem \ref{p<2.1} that provides  information on complemented subspaces of $L_p$ spaces for all $1<p<\infty$.  The proof is similar to the proof of  $(4)\implies (1)$ in Theorem \ref{p<2.1} but uses the well known consequence of  Lindenstrauss' \cite{lin} work that every reflexive space $X$ has an $M$-basis; that is, a biorthogonal system $\{x_\alpha, x_\alpha^*\}_{\alpha\in A}$ such that the linear span of $\{x_\alpha\}_{\alpha\in A}$ is dense in $X$ and $\{x_\alpha^*\}_{\alpha\in A}$ separates points of $X$ (which, in view of the reflexivity of $X$, is the same as saying that the linear span of $\{x_\alpha^*\}_{\alpha\in A}$ is dense in $X^*$); see, for example, \cite[Theorem 4.2]{ziz}.

  In the proof of Theorem \ref{p<2.1+}  we use the following simple lemma, which is a special case of a result well known to people who work in non separable Banach space theory.
  To avoid cumbersome notation in the statement and proof of the lemma, we introduce here some notation.  Fix $1<p<\infty$  and consider a space $L_p = (\sum_{\beta\in B} L_p\{-1,1\}^{\aleph_\beta})_p$. We call a projection $P$ on $L_p$ of the form $\sum_{\beta\in B'}   \Bbb{E}_{C_\b}$ with $B'\subset B$ and $C_\beta$ a non empty subset of $\aleph_\beta$ for $\beta\in B'$ a {\sl standard} projection. Here $\Bbb{E}_{C_\b}$ is the conditional expectation projection from  $L_p\{-1,1\}^{\aleph_\beta}$ onto $L_p\{-1,1\}^{C_\beta}$.
  So if $P$ is a standard projection on $L_p$, then $P^*$ is the obvious standard projection on $L_{p^*}$.  Given standard projections $P$ and $Q$ on $L_p$, write $P\le Q$ if $QP=P$. So if $P\le Q$, then $QP=PQ=P$.  Evidently if  $P=\sum_{\beta\in B'}  \Bbb{E}_{C_\beta}$ and $Q= \sum_{\beta\in B''}  \Bbb{E}_{C'_\beta}$, then $P\le Q$ iff $B'\subset B''$ and for all $\beta\in B'$ we have $C_\beta \subset C'_\beta$. Notice that  the density character of the range of a standard projection $P=\sum_{\beta\in B}    \Bbb{E}_{C_\beta}$  is $\aleph_0 + \sum_{\beta\in B} |C_\beta|$.

  \begin{lm}\label{p<2.1+lemma} Let $\{x_\alpha\}_{\alpha\in A}$ be an $M$-basis for a subspace $X$ of \hfill\break $L_p:=  (\sum_{\beta\in B} L_p\{-1,1\}^{\aleph_\beta})_p$ and let $P$ be a standard projection on $L_p$.  Then there is a standard projection $Q\ge P$ on $L_p$ so that the density character $\aleph$ of the range of $Q$ is the same as the density character of the range of $P$ and so that there is a subset $A'\subset A$ with $|A'| \le \aleph$ and $Qx_\a =x_\a$ for $\a \in A'$ and $Qx_\a =0$ for $\a \not\in A'$.
\end{lm}

 \pf The main point is the obvious fact that if $T$ is an operator from $X$ into some Banach space and the range of $T^*$ is contained in the closed  span of $\{x^*_\a\}_{a\in A'}$, then for all $\a \not\in A'$ we have $Tx_\a=0$.  Construct by induction an increasing sequence $P_n$ of standard projections on $L_p$ and subsets \hfill\break  $A_1\subset A_2\subset \dots$ of $A$ so that for each $n$ we have
 \begin{itemize}
 \item The range of $({P_n}_{|X})^*$ is contained in the closed span of $\{x^*_\a\}_{a\in A_n}$.
 \item $P_{n+1}x_\a = x_\a$ for all $\a\in A_n$.
 \end{itemize}
 Use the ``main point"  to check that the sequence $P_n$ converges strongly to a standard projection $Q$  that satisfies the conclusions of the lemma.  \qed

 \begin{thm}\label{p<2.1+} Let $X$ be a subspace of $L_p=  (\sum_{\beta\in B} L_p\{-1,1\}^{\aleph_\beta})_p$, $1<p<2$,  and let $\aleph$ be an uncountable cardinal.  If $\ell_p(\aleph)$ is not isomorphic to a subspace of $X$, then for every $\e>0$ there is a subset $B'$ of $B$ with $|B'|<\aleph$ so that the restriction projection $P_{B'}$ of $L_p$ onto $ (\sum_{\beta\in B'} L_p\{-1,1\}^{\aleph_\beta})_p$ satisfies $\|P_{B'}x -x\| \le \e \|x\|$ for every $x\in X$.
 \end{thm}

 \pf Since $X$ has an $M$-basis we deduce from Lemma \ref{p<2.1+lemma} that if $P$ is a standard projection on $L_p$, then there is a standard projection $Q$ on $L_p$ so that $Q\ge P$ and $QX\subset X$ and the range of $Q$ has the same density character as the range of $P$. We now basically just repeat the proof of $(4)\implies (1)$ in Theorem \ref{p<2.1}, but input this consequence of Lemma \ref{p<2.1+lemma}.

Assume that the conclusion of Theorem \ref{p<2.1+} is false for a certain $\e>0$.
Take a set $M:= \{x_\g : \g\in \G\}$ of unit vectors in $X$ maximal with respect to the properties
\begin{itemize}
\item[(i)] The collection $M$ is a GMD set.
\item[(ii)] There are disjoint countable subsets $B_\g$ of $B$ such that for all $\g$ in $\G$, $\|P_{B_\g} x_\g \| > \e/2$.
\end{itemize}
To get   a contradiction, just as in the proofs of Theorem \ref{p<2} and Theorem \ref{p<2.1},  it is enough to verify that $|\G|\ge \aleph$. Assume instead that $\aleph > |\G|$. We might as well assume, by enlarging the $B_\g$'s, that for $\b\not\in B' := \cup_{\g\in\G} B_\g$, we have  $P_\b x_\g = 0 $ for all $\g$ in $\G$.  The set $B'$ could be replaced with any superset that has cardinality at most $|\G|$.

For each $\b$ in $B'$ and $\g$ in $\G$, let $C_{\b,\g}$ be the (countable) set of coordinates of $\{-1,1\}^{\aleph_\b}$ on which $P_\b x_\g$ depends, and let $C_\b := \cup_{\g\in\G} C_{\b,\g}$.  Then $|C_\b|\le |\G|$ and hence the density character of $ L_p\{-1,1\}^{C_\b}$ is at most $|\G|$.  We could as well replace $C_\b$ with any subset of $\aleph_\b$ of cardinality at most $|\G|$ that contains this $C_\b$.
 Let   $\Bbb{E}_\b$  be the (norm one) conditional expectation projection from $L_p\{-1,1\}^{\aleph_\b}$  onto $L_p\{-1,1\}^{C_\b}$ and let $P:= (\sum_{\b\in B'} \Bbb{E}_\b)_p$ be the direct sum of these projections; so $P$ is a standard projection on $L_p$. Note that the density character of the range of  $P$ is at most $|\G|<\aleph$. By the consequence of Lemma  \ref{p<2.1+lemma} we mentioned above, there is a standard projection $Q\ge P$ so that the density characters of the ranges of $Q$ and $P$ are the same and $QX\subset X$.  By the replacement comments made above, without loss of generality we can  avoid
  introducing additional notation and just assume that $P=Q$.

 We are now ready for the punch line.  If $y$ is a unit vector in the subspace  $(I-P)X$ of $X$, then $M\cup \{y\}$ is still a GMD set.  Consequently, by maximality of $M$ we have $\|(I-P_{B'})y\| \le \e/2$.  Thus if $x$ is a unit vector in $X$ we have
 since $P\le P_{B'}$ that
 $$
 \|x-P_{B'}x\| = \|(I-P_{B'})(I-P)x\|\le (\e/2)\|(I-P)x\|\le \e\|x\|.  \ \ \  \qed
 $$

Our main reason for proving Theorem \ref{p<2.1+} is that it gives as a corollary a verification of the conjecture mentioned in the Introduction for complemented subspaces of $L_p$ spaces, $2<p<\infty$.

\begin{co}\label{complemented}  Assume that $X$ is a complemented subspace of  some $L_p$ space (say, with a projection of norm $\lambda$). Let $\aleph$ be an infinite cardinal and assume that $\ell_p(\aleph)$ is not isomorphic to a subspace of $X$.  Then
\begin{itemize}
\item[(1)] If $1<p<2$, then  for all $\e>0$ there is $\aleph'<\aleph$ such that  the space $X$ is $(1+\e)$-isomorphic to a $(1+\e)\lambda$-complemented subspace of some $L_p(\aleph')$ space.
\item[(2)] If $2<p<\infty$, then there is $\aleph'<\aleph$ such that  $X$ is isomorphic to a  complemented subspace of some $L_p(\aleph')$ space.
\end{itemize}
\end{co}

 \pf In view of Maharam's theorem,  (1) is immediate from Theorem \ref{p<2.1+}  and the principle of small perturbations.  Conclusion (2) follows from (1) by duality and the equivalence of  (2)
 and (3) in Theorem \ref{p<2.1}.  \qed

 \begin{tabular}{l}
W.~B.~Johnson\\
Department of Mathematics\\
Texas A\&M University\\
College Station, TX  77843 U.S.A.\\
{\tt johnson@math.tamu.edu}
\\
\end{tabular}

\bigskip

\begin{tabular}{l}
G.~Schechtman\\
Department of Mathematics\\
Weizmann Institute of Science\\
Rehovot, Israel\\
{\tt gideon@weizmann.ac.il}\\
\end{tabular}

\end{document}